\documentclass[12pt,reqno]{amsart}

\usepackage{amssymb}
\usepackage{amsmath}
\usepackage{amsthm}
\usepackage{mathtools}
\usepackage{eucal}
\usepackage{color}
\usepackage{amsfonts}

\newcommand{\R}{\mathbb R}

\newcommand{\p}{\partial}

\newcommand{\sgn}{\text{sgn}}

\newcommand{\norm}[1]{\left\| #1 \right\|}

\newcommand{\jd}{\rangle}

\newcommand{\ji}{\langle}

\numberwithin{equation}{section}
\newcommand{\be}{\begin{equation}}
\newcommand{\ee}{\end{equation}}
\newcommand{\bp}{\begin{proof}}
\newcommand{\ep}{\end{proof}}
\newcommand{\bel}{\begin{equation}\label}
\newcommand{\eeq}{\end{equation}}
\newtheorem{thm}{Theorem}[section]

\newtheorem{lem}[thm]{Lemma}

\theoremstyle{remark}
\newtheorem{rem}{Remark}[section]

\numberwithin{equation}{section}

\begin{document}

\title{On local energy decay for solutions to the b-family of peakon equations}
\author{C. Hong}
\address[C. Hong]{Department  of Mathematics\\
University of California\\
Santa Barbara, CA 93106\\
USA.}
\email{christianhong@ucsb.edu}


\begin{abstract} 
    This work is concerned with the long time behavior of solutions to the $b$-family of peakon equations. We prove local energy decay of global solutions under suitable hypotheses. Assuming the global bound of the $H^1(\R)$ norm, we show local energy decay along sequences of time in an expanding region around the origin. If we assume nonnegativity for an initial momentum, $m_0(x) = (u - \p_x^2 u)(0,x)$, we show strict decay on a similar region. Moreover, we show decay in an exterior region given the same nonnegativity condition.
\end{abstract}

\maketitle

\section{Introduction}

In this work, we consider the Camassa-Holm (CH) equation posed on the real line

\be
\label{loc-CH}
\p_t u - \p_{xxt}^3 u + 3u \p_x u - u \p_x^3 u - 2 \p_x u \p_x^2 u =0, \, \, (t, x) \in \R^2.
\ee

The CH equation was first noticed by Fuchssteiner and Fokas \cite{FF1981} when they studied hereditary symmetries. Later it was derived by Camassa and Holm \cite{CH1993} as an alternative integrable model for shallow water waves in a channel. This model was also deduced by Dai \cite{D1998} in the context of elastic rods. Among the notable properties of this system are the existence of a bi-Hamiltonian structure and the existence of peaked traveling wave solutions called peakons.

As discussed in \cite{CH1993}, the CH equation is completely integrable and posseses a Lax pair, which gives infinitely many conservation laws. The first two laws are given by

\be
\label{conserved}
I(u) = \int u \, dx, \, \, \,  E(u) = \int \left( u^2 + (\p_x u)^2 \right) \, dx.
\ee

Li and Oliver \cite{LO2000} and Rodriguez \cite{R2001} showed that the initial value problem (IVP) associated to the CH equationis (strongly) locally well-posed in $H^2(\R)$, $s > 3/2$. The second conservation law hints at the existence of global solutions in $H^1(\R)$ in the weak sense. Constantin and Molinet \cite{CM2000} managed to prove one can get such solutions that are continuous in time with respect to the $H^1(\R)$ norm under a nonnegativity condition on the quantiy $m = u - \p_x^2 u$. McKean \cite{Mc2004} then showed that solution blow-up occurs if and only if $m$ has a positive part to the left of some negative part of itself. In addition, one can establish a type of local well posedness for solutions in $H^1(\R) \cap W^{1, \infty}(\R)$, which is a class that includes our peakon solutions. This was proven by Linares, Ponce, and Sideris in \cite{LPS2016}.

Constantin \cite{C2000} showed that classical solutions will only blow-up as wave-breaking, which means that the slope of a solution may diverge at a point. Moreover, Bressan and Constantin \cite{BC2007} managed to construct global conservative solutions in $H^1(\R)$ for any initial condition. Such solutions can experience wave breakdown, but can be extended globally in time after they blow up.

Using that $(1 - \p_x^2)^{-1}g = f$ implies $Gf = \frac 1 2 e^{-|x|} * f = g$, one gets the conservative nonlocal form of the CH equation.

\be
\label{CH}
\p_t u + \p_x \left(\frac 1 2 u^2 + G \left( u^2 + \frac 1 2 (\p_x u)^2 \right) \right)=0, \, \, (t,x) \in \R \times \R.
\ee

This weaker formulation justifies the validity of peakon solutions, which are peaked traveling wave equations to our equation. These peakons take the general form

\be
u(t,x) = cQ(x -ct), \, \, Q(x) = c e^{-|x|}, \, \, c \in \R.
\ee

Notice that $Q \notin H^{3/2}(\R)$.

As it was mentioned before, solutions to the CH equation may experience blow-up in finite time via wave breaking and still can be extended to $H^1(\R)$ global solutions. Such behavior does not occur for similar models. For example, some solutions to the mass critical gKdV equation can blow up in $H^1(\R)$ by having mass escape to infinity in finite time. This is fully described by Martel and Merle \cite{MM2000}. Other examples, such as the KdV equation, possess global solutions in $H^1(\R)$. 

The CH equation demonstrates elastic collisions between peakons, much like soliton behavior for the KdV equation due to their integrability. The stability of single peakons was established by Constantin and Strauss \cite{CS2000}. As hypothesized during its derivation, asymptotic stability of trains of peakons solutions have been established by Molinet in \cite{M2018}, which aligns with asymptotic stability of solitons for the KdV equation. The speed of these peakons are determined by their amplitude, which is bounded by their $H^1(\R)$ by the Sobolev embedding. This implies that solutions should decay in the $H^1(\R)$ norm in the shrinking region $|x| \gtrsim \beta t$ for some $\beta > 0$.

A train of peakons should decay around the origin since they will carry their mass away from the origin at constant speed. This combined with asymptotic stability results implies decay on some sort of expanding interval centered at the origin. Alejo, Cortez, Kwak, and Mu\~{n}oz \cite{ACKM2019} proved this inside the interval approximately of the size $|x| \lesssim t^{1/2}$. The relevant theorem is

\begin{thm}
    \label{Decay}
    \cite{ACKM2019} Let $u$ be a nontrivial global solution to \eqref{CH}, such that $u \in C([0, \infty]:H^1(\R)) \cap L^\infty (\R:L^1(\R))$. Then we have

    \be
    \lim_{t \to +\infty} \| u(t) \|_{H^1(\Lambda_{1/2}(t))} = 0
    \ee

    where
    \be
    \label{Ib}
    \Lambda_c(t) \coloneqq \left( -\frac{|t|^c} {\log |t|}, \frac {|t|^c} {\log |t|} \right).
    \ee

\end{thm}

\begin{rem}
Note that \eqref{conserved} and the global existence of solutions with nonnegative initial data \cite{C2000} imply that the hypotheses, $u \in C([0, \infty]:H^1(\R)) \cap L^\infty (\R:L^1(\R))$, are satisfied for $m_0(x) \ge 0$, and $m_0 \in L^1(\R)$.
\end{rem}

\begin{rem}
The above theorem was also proved for the $L^2$-norm for the Degasperis-Procesi (DP) equation, which was first discovered by Degasperis and Gaeta \cite{DG1999}. This has the form

\be
\label{DP}
\p_t u + \p_x\left( \frac 1 2 u^2 + \frac 3 2 G u^2 \right) = 0.
\ee

This equation also has the same peakon solutions, and is a bi-Hamiltonian system. In addition, \eqref{DP} allows for shock peakon solutions, which are like peakons except they are discontinuous at a point. This behavior is possible since weak solutions to \eqref{DP} may simply be $L^2(\R)$ instead of $H^1(\R)$. These solutions take the form

\be
\label{shocks}
u(x,t) = \frac 1 {t + k} \sgn (x) e^{-|x|}, \, \, k > 0.
\ee
\end{rem}

In \cite{ACKM2019}, one finds a similar result by accepting a growth in time of the $L^1(\R)$ norm of solutions, as opposed to requiring a global bound. In particular, one allows polynomial growth. The result is

\begin{thm}
    \label{WeakDecay}
    Let $u \in C([0, \infty): H^1(\R))$ be a solution to \eqref{CH} or \eqref{DP} such that

    \be
    \int |u(t, x)| \, dx \lesssim \ji t \jd^a.
    \ee

    For $0 \le a < 1$, and $c \coloneqq 1 - a$, we have
    \begin{enumerate}
        \item For solutions to CH \eqref{CH}
              \be
              \liminf_{t \to \infty} \| u(t) \|_{H^1(\Lambda_c(t))} = 0.
              \ee

        \item For solutions to DP \eqref{DP}
              \be
              \liminf_{t \to \infty} \| u(t) \|_{L^2(\Lambda_c(t))} = 0.
              \ee

    \end{enumerate}

\end{thm}

If $a < 1/2$, this result holds for a larger expanding interval, but it only gets a $\liminf$ result. The proof extends to the $b$-family of equations, which are discussed below.

The CH equation and the DP equation belong to the $b$-family of equations

\be
\p_t u - \p_{xxt}^3 u + (b+1) u \p_x u - b \p_x u \p_x^2 u - u \p_x^3 u = 0,
\ee

or equivalently

\be
\label{bfam}
    \p_t u + \p_x \left( \frac 1 2 u^2 + G \left(\frac b 2 u^2 + \frac {3-b} 2 (\p_x u )
    ^2 \right) \right) = 0, \, b \in (0, 3),
\ee

The CH equation is represented when $b=2$ while the DP equation is represented when $b=3$. Note that all of these equations admit peakon solutions, but only the DP equation admits breakdown solutions. A caveat in this relation is that only the CH equation and the DP equation are integrable. A question we seek to answer in this paper is whether results similar to those in \cite{ACKM2019} apply to \eqref{bfam} solutions in general since the integrability of these equations are not directly applied in the proof of Theorem \ref{Decay}. These hypotheses are not unfounded since it was shown in \cite{Ch2023} that peakon solutions are stable for $b > 1$. The same paper found numerical evidence of the formation of traveling wave-like behavior for all $0 < b < 3$, although the system behaved more like the viscous Burger's equation when $b < 1$. This also begs the question whether we will have decay in the region $|x| \gtrsim \beta t$. Since \eqref{bfam} is not integrable for $b \ne 2,3$, one does not have access to inverse scattering techniques to investigate this behavior. In addition, we are also interested in a similar result that does not require any assumptions on the longtime behavior of the $L^1(\R)$ norm of the solution. This is a natural question since there exists global solutions to the CH equation that are bounded in $H^1(\R)$, but have varying behavior in the $L^1(\R)$ norm. This of course extends to the b-family.

Since the $b$-family equations are not integrable in general, we have to take advantage of different conserved quantities than what were used in the theorems mentioned above. Consider the quantity $m(t,x) = (1- \p_x^2)u(t,x)$. We can rewrite \eqref{bfam} in terms of $m$. As stated in \cite{GLT2008}, This gives us

\be
\label{cons-bfam}
\p_t m + u \p_x m + b m\p_x u = 0.
\ee

It can be shown that the sign of $m(t, x)$ is determined by the sign of $m(0, x)$. In particular, if $m(0, x)$ is nonnegative, $m(t,x)$ remains nonnegative. Using the Green's function of $G$, we learn that $m(0, x) = m_0(x) \ge 0$ implies $u(t, x) \ge 0$ since $\frac 1 2 e^{-|x|}$ is nonnegative.

Like the $b=2$ CH case, $m$ is an important quantity for the analysis of these equations. For example, it was shown that $u \in C([0, \infty) :H^2(\R))$ is a global solution to \eqref{CH} if $m_0$ is signed. This is similar to the behavior of \eqref{loc-CH}. Breakdown of solutions even occur in the same way, with the slope of our solutions diverging. The proof of these results can be found in \cite{GLT2008}. Moreover, $\int m\, dx$ is a conserved quantity for solutions of \eqref{bfam}. This means that for all $u$ that solves \eqref{bfam}, we have

\be
\frac d {dt} H(t) = \frac d {dt} \int_{-\infty}^\infty m \, dx = \frac d {dt} \int_{-\infty}^\infty (u - \p_x^2 u)(x,t) \, dx = 0.
\ee

Note that \eqref{CH} and \eqref{DP} have other integrals of motion that are very useful. In particular, solutions to the CH equation \eqref{CH} and the DP equation \eqref{DP} have constants of motion that are equivalent to the $H^1(\R)$ norm. Since the $b$-family lacks access to these conservation laws, it is difficult to obtain the exact same results for $b \ne 2, 3$. This is in spite of the fact that \cite{Ch2023} found very similar behavior numerically when $b > 1$. The techniques used in this paper seem unable to take advantage of this fact, so it is possible that further results can be found by taking advantage of whatever properties that are tied to $b > 1$.

Next, we will discuss our results.

If one weakens their assumptions on regularity of their solutions but without any hypotheses on the asymptotic behavior of the $L^1(\R)$ norm of the solution, one can still get a weak decay result for solutions. In particular, one can get $\liminf_{t \to \infty} \| u(t) \|_{H^1(\Lambda_{2/3}(t))} = 0$ on expanding inervals. 
\begin{thm}
    \label{TH3}
    Let $u \in C([0, \infty) : H^1(\R)) \cap L^\infty([0, \infty) : H^1(\R))$ be a global solution to \eqref{bfam}. Then we have

    \be
    \liminf_{t \to \infty} \| u(t) \|_{H^1(\Lambda_{2/3}(t))}.
    \ee

\end{thm}

 If we instead restrict the $L^1$ norm of $m_0(x)$, we get a similar result to Theorem \ref{Decay}. This result is

\begin{thm}
    \label{TH1}
    Let $u \in C([0, \infty) : H^2(\R))$ be a global solution to \eqref{bfam} such that $m_0 \ge 0$, and $m_0 \in L^1(\R)$, then

    \be
    \lim_{t \to \infty} \| m(t) \|_{L^1(\Lambda_{1/2}(t))} = \lim_{t \to \infty} \| u(t) \|_{L^1(\Lambda_{1/2}(t))} = L
    \ee

    for some $L \in \R$.

\end{thm}

\begin{rem}
    There is a similar decay result for the Zakharov-Kuznetsov equation that was proven by Mendez, Mu\~{n} oz, Poblete, and Pozo \cite{MMPP2021}. Our argument was inspired by their approach.
\end{rem}

\begin{rem}
Note that we require $u$ to be in $H^2(\R)$ for all time. This is because we have to make use of $m = u - u_{xx}$ to gain control over the $H^1(\R)$ norm for $u$. For the CH equation, this is not necessary because $\| u \|_{H^1}$ is a conserved quantity for the CH equation. Moreover, we do not know what the limit $L$ will be.
\end{rem}

Finally, we can show that signed solutions to \eqref{bfam} decay to zero on $(\sigma t, \infty)$ and $(-\infty, -\sigma t)$. In essence, this implies the traveling waves we expect from these PDEs to have a velocity capped by some quantity related to the initial data.

\begin{thm}
    \label{TH4}
    Let $u \in C([0, \infty) : H^2(\R))$ be a global solution to \eqref{bfam} such that $m_0 \ge 0$, and $m_0 \in L^1(\R)$. Then there exists $\sigma > 0$ such that 
    \[ \lim_{t \to \infty} \| u \|_{W^{1,p}(\overline\lambda(t))} \to 0, \]
    where
    \[ \overline\lambda(t) = (\sigma t, \infty). \]
    One can show a similar result and proof for $\tilde\lambda(t) = (-\infty, -\sigma t)$ when $m_0 \le 0$.

\end{thm}

We once again require $m$ to be signed in order to easily bound the $L^1(\R)$ norm of $m$. We follow a similar proof for the gBBM equation shown in \cite{K2019}.

\section{Preliminary Estimates}

We need to bound the $H^1(\R)$ norm on $u$ in time for solutions to \eqref{bfam}. However, with the exception of $b=2$, the $H^1(\R)$ norm is not a conserved quantity. We can take advantage of $m = u - u_{xx}$ and $G$ to get a bound on the $H^1(\R)$ norm instead. To do this, note that 

\[u = G m = \frac 1 2 e^{-|x|} * m, \, \, \, \p_x u = \p_x G m = - \frac {\sgn(x)} 2 e^{-|x|} * m. \]

One speical case of Young's convolution inequality tells us that for $1 \le p \le \infty$, $f \in L^1$, and $g \in L^p$.
\[ \| f * g \|_p \le \| f \|_1 \| g \|_p. \]

We have $\|\frac 1 2 e^{-|x|}\|_p = \| - \frac {\sgn(x)} 2 e^{-|x|} \|_p \le 1$. This then implies

\be
\label{p-estimate}
 \| u \|_p \le \| m \|_1, \, \, \, \| \p_x u \|_p \le \|m \|_1, \, 1 \le p \le \infty.
\ee

Recall that $\int m \, dx$ is a conserved quantity, and $m$ remains signed if $m(0, x)$ is signed. This gives a condition to bound the $L^1(\R)$ norm of $m$, which can then be used to bound the $H^1(\R)$ norm of $u$ for all time. Moreover, a decay on $\| m \|_1$ is equal to a decay in $\| u \|_{H^1(\R)}$.

\subsection{Estimates on $\lambda(t)$ and $G$}

We will be taking advantage of a couple of functions. Let $b \in [0, 1)$. We will define the following weight function $\phi(x)$.

\be
\label{weight}
\phi(x) = \begin{cases} e^x -1 & x \le 0 \\ 1 - e^{-x} & x > 0
\end{cases}, \, \phi'(x) = e^{-|x|}, \, \phi''(x) = -\sgn(x) e^{-|x|}.
\ee

We will also use the following functions $\lambda_c(t), \mu_c(t)$. $\lambda_c(t)$ is related to the growing interval where we measure decay, while $\mu_c(t)$ is necessary to prove Theorem \ref{TH3}.

\be
\label{lu}
\lambda_c(t) = \frac{t^c}{\log t}, \quad \mu_c(t) = t^{1-c}\log^2 t, \quad \mu_c(t) \lambda_c(t) = t \log t.
\ee

In the case $c=0$, we adopt the convention $\lambda_c(t) = 1$. When $c > 0$, we have

\be
\label{lu-deriv}
\lambda_c'(t) = \frac{t^{c-1}}{\log t} \left(c - \frac 1 {\log t} \right), \quad \frac{\lambda_c'(t)}{\lambda_c(t)} = \frac 1 t \left(c - \frac 1 {\log t} \right), \quad \frac {\mu_c'(t)}{\mu_c(t)} = \frac {1-c} t + \frac {2 \log t} t.
\ee

We will also take advantage of the fact that $(1 - \p_x^2)^{-1}f = Gf = \frac 1 2 e^{-|x|} * f$. Since $\frac 1 2 e^{-|x|} > 0$, we get

\be
f \ge g \implies Gf \ge Gg, \quad f \ge 0 \implies Gf \ge 0.
\ee

We will also need the following Lemma
\begin{lem}
    \label{exp-est}

    \be
    \begin{aligned}
        G\phi'(x/\lambda) = \frac{\lambda^2}{\lambda^2-1}e^{-|x|/\lambda} - \frac{\lambda}{\lambda^2-1}e^{-|x|} \approx \phi'(x/\lambda)
    \end{aligned}
    \ee
    when $\lambda \gg 1$.

\end{lem}
\bp
Let's consider when $x \le 0$
\be
\begin{aligned}
    2G\phi'(x/\lambda) & = \int_{-\infty}^\infty e^{-|x-y|}e^{-|y|/\lambda} \, dy                                                                                       \\
                                           & = \int_{-\infty}^{x} e^{y-x} e^{y/\lambda}\, dy + \int_x^0 e^{x-y}e^{y/\lambda} \, dy + \int_0^\infty e^{x-y}e^{-y/\lambda} \, dy              \\
                                           & = e^{-x} \int_{-\infty}^{x} e^{y+y/\lambda}\, dy + e^x\int_x^0 e^{-y+y/\lambda} \, dy + e^x\int_0^\infty e^{-y-y/\lambda} \, dy                \\
                                           & = e^{-x} \left. \frac 1 {1 + 1/\lambda} e^{y+y/\lambda} \right|_{-\infty}^x + e^x \left. \frac 1 { 1/\lambda-1} e^{-y+y/\lambda} \right|_{x}^0 \\
                                           & + e^x \left. \frac 1 {-1 - 1/\lambda} e^{-y-y/\lambda} \right|_0^\infty                                                                        \\
                                           & = \frac 1 {1 + 1/\lambda} e^{x/\lambda} + \frac 1 { 1/\lambda-1} (e^{x} - e^{x/\lambda}) + \frac 1 {1 + 1/\lambda} e^x                         \\
                                           & = \frac{2\lambda^2}{\lambda^2-1} e^{x/\lambda} - \frac{2\lambda}{\lambda^2-1}e^x
\end{aligned}
\ee
We get a similar computation for $x > 0$. Note that since $\lambda > 1$, $e^{x/\lambda}$ dominates $e^x$, so we get the result. Moreover, as $\lambda \to \infty$, $G\phi'(x/\lambda) \approx \frac {\lambda^2} {\lambda^2-1} \to \phi'(x /\lambda)$ uniformly.
\ep

\section{Proof of Theorem \ref{TH3}}

\subsection{Functional Estimate}
Consider the functional

\be
\mathcal M_q(t) = \frac 1 {\mu_c(t)} \int \phi \left( \frac x {\lambda_c(t)} \right) \phi' \left( \frac x {\lambda_c^q(t)} \right) u(t, x) \, dx
\ee

where $0 < c < 1$, and $q > 1$. We will find a tighter bound on $b$ as we prove the upcoming lemma. To use this functional, we will need to show that it is bounded above after sufficient time. Let us assume that $u \in C([0, \infty): H^1(\R))$. Since $\phi \in L^\infty(R)$, and $\phi' \in H^1(\R)$, we know that $\mathcal M_q(t) < +\infty$ for all $t > 2$. In fact,

\be
\begin{aligned}
    | \mathcal M_q(t)| & \le \frac 1 {\mu_c(t)} \norm{\phi \left( \frac \cdot {\lambda_c(t)} \right)}_{L^\infty} \norm{\phi' \left( \frac \cdot {\lambda_c^q(t)} \right)}_{L^2} \norm{u(t)}_{L^2} \\
              & = \frac {\lambda^{q/2}(t)} {\mu_c(t)} \norm{\phi }_{L^\infty} \norm{\phi' }_{L^2} \norm{u(t)}_{L^2}                                                                  \\
              & = C \frac {t^{cq/2}} {t^{1-c} \log^{2} t \log^{q/2}}\norm{u(t)}_{L^2} = C \frac 1 {t^{(2-2b-bq)/2}\log^{(4+q)/2} t}\norm{u(t)}_{L^2}.
\end{aligned}
\ee

$\norm{u(t)}_{L^2} \le \norm{u(t)}_{H^1}$, which is bounded by hypothesis, so if we let $b \le \frac 2 {2+q}$, we get

\be
\sup_{t \gg 1}|\mathcal M_q(t)| = M < +\infty.
\ee

We will use this fact in conjunction with the derivative of $\mathcal M_q (t)$. If we take the derivative of $\mathcal M_q(t)$ with respect to time, (1.1) and integrating by parts tells us that

\be
\label{M-deriv}
\begin{aligned}
    \frac d {dt} \mathcal M_q (t) & = -\frac {\mu_c'(t)} {\mu_c^2(t)} \mathcal M_q(t) - \frac 1 {\mu_c(t)} \frac {\lambda_c'(t)} {\lambda_c(t)} \int \frac x {\lambda_c(t)} \phi' \left( \frac x {\lambda_c(t)} \right) \phi' \left( \frac x {\lambda_c^q(t)} \right) u(t, x) \, dx     \\
                                  & - \frac q {\mu_c(t)} \frac {\lambda_c'(t)} {\lambda_c(t)} \int\ \frac x {\lambda_c^q(t)} \phi \left( \frac x {\lambda_c(t)} \right) \phi'' \left( \frac x {\lambda_c^q(t)} \right) u(t, x) \, dx                                              \\
                                  & + \frac 1 {\mu_c(t)\lambda_c(t)} \int \phi' \left( \frac x {\lambda_c(t)} \right) \phi' \left( \frac x {\lambda_c^q(t)} \right) \left[\frac {u^2} 2 + G \left( \frac b 2 u^2 + \frac {3-b} 2 (\p_x u)^2\right)\right] \, dx   \\
                                  & + \frac 1 {\mu_c(t)\lambda_c^q(t)}  \int\phi \left( \frac x {\lambda_c(t)} \right) \phi'' \left( \frac x {\lambda_c^q(t)} \right) \left[\frac {u^2} 2 + G \left( \frac b 2 u^2 + \frac {3-b} 2 (\p_x u)^2\right)\right] \, dx \\
                                  & = -M_1 -M_2 - M_3 + M_4 + M_5.
\end{aligned}
\ee

Now we can prove the following lemma

\begin{lem}
    Let $u$ be a nontrivial global solution of \eqref{bfam}, such that $u \in C([0, \infty): H^1(\R)) \cap L^\infty([0, \infty): H^1(\R))$. If $b < 2/3$, we get
    \be
    \int_{10}^\infty \frac 1 {t \log t} \int \phi'\left(\frac x {\lambda_c(t)} \right) \phi'\left(\frac x {\lambda_c^q(t)} \right) (u^2 + (\p_x u)^2) \, dxdt < +\infty.
    \ee
\end{lem}

\bp
Note that $M_1 = \frac {\mu_c'(t)} {\mu_c^2(t)} \mathcal M_q(t)$ is integrable for $t \gg 1$. This is because $\frac {\mu_c'(t)} {\mu_c^2(t)} = \frac {1 + \log t}{t^2 \log^2 t}$, which is integrable for $t \gg 1$, and $\mathcal M_q(t)$ is bounded for $t \gg 1$ as well. Similar reasoning tells us that

\be
\begin{aligned}
    M_5 & = \frac 1 {t^{1-c+cq}\log^{2-q}t}  \int\phi \left( \frac x {\lambda_c(t)} \right) \phi'' \left( \frac x {\lambda_c^q(t)} \right) \left[\frac {u^2} 2 + G \left( \frac b 2 u^2 + \frac {3-b} 2 (\p_x u)^2\right)\right] \, dx \\
        & \le C \frac 1 {t^{1-c+cq}\log^{2-q}t} \norm{\phi\phi''}_{L^\infty} \norm{u}_{H^1}^2.
\end{aligned}
\ee

The last inequality follows from the fact that $Gf = \frac 1 2e^{-|x|} *f$. $M_2$ is then integrable because $q > 1$.\\
Young's inequality implies that

\be
\label{M2}
\begin{aligned}
    \left|M_2\right| & \le \left|\frac{\lambda_c'(t)}{2\epsilon\mu_c(t)} \right|\int \left(\frac x {\lambda_c(t)} \right)^2 \left[ \phi' \left( \frac x {\lambda_c(t)} \right)\phi' \left( \frac x {\lambda_c^q(t)} \right)\right] \, \frac {dx}{\lambda_c(t)} \\&+ \left|\frac{ \epsilon\lambda_c'(t)}{2\mu_c(t)\lambda_c(t)}\right| \int \left[ \phi' \left( \frac x {\lambda_c(t)} \right)\phi' \left( \frac x {\lambda_c^q(t)} \right)\right] u^2 \, dx \\
                     & = \left|\frac{(\lambda_c'(t))^2}{2\mu_c(t)} \right|\int \left(\frac x {\lambda_c(t)} \right)^2 \left[ \phi' \left( \frac x {\lambda_c(t)} \right)\phi' \left( \frac x {\lambda_c^q(t)} \right)\right] \, \frac {dx}{\lambda_c(t)}       \\&+ \left|\frac{ 1}{2\mu_c(t)\lambda_c(t)}\right| \int \left[ \phi' \left( \frac x {\lambda_c(t)} \right)\phi' \left( \frac x {\lambda_c^q(t)} \right)\right] u^2 \, dx
\end{aligned}
\ee

when $\epsilon = \frac 1 {\lambda_c'(t)}$.\\
Similar reasoning tells us that

\be
\label{M3}
\begin{aligned}
    \left|M_3\right| & \le \left|\frac q {2\mu_c(t)} \frac {\lambda_c'(t)} {\epsilon\lambda_c(t)} \right| \int \left(\frac x {\lambda_c^q(t)} \right)^2 \left|\phi \left( \frac x {\lambda_c(t)} \right)\right| \left[\phi'' \left( \frac x {\lambda_c^q(t)} \right)\right]^2  \, dx \\ &+ \left|\frac {\epsilon q} {2\mu_c(t)} \frac {\lambda_c'(t)} {\lambda_c(t)} \right|\int  \left|\phi \left( \frac x {\lambda_c(t)} \right)\right|  u^2 \, dx \\
                     & = \left|\frac q {2\mu_c(t)} \frac {\lambda_c'(t)} {t^{1-c}\lambda_c(t)} \right| \int \left(\frac x {\lambda_c^q(t)} \right)^2 \left|\phi \left( \frac x {\lambda_c(t)} \right)\right|\left[ \phi'' \left( \frac x {\lambda_c^q(t)} \right)\right]^2  \, dx    \\ &+ \left|\frac {qt^{1-c}} {2\mu_c(t)} \frac {\lambda_c'(t)} {\lambda_c(t)} \right|\int  \left|\phi \left( \frac x {\lambda_c(t)} \right)\right|  u^2 \, dx
\end{aligned}
\ee

when $\epsilon = t^{1-c}$.\\
In \eqref{M2}, note that $|\phi'(x)| \le 1$. We can then use \eqref{weight}, \eqref{lu-deriv}, and the change of variables $y=\frac x {\lambda_c(t)}$ to get

Say 

\be
M_3^{(1)} = \frac{(\lambda_c'(t))^2}{2\mu_c(t)} \int \left(\frac x {\lambda_c(t)} \right)^2 \left[ \phi' \left( \frac x {\lambda_c(t)} \right)\phi' \left( \frac x {\lambda_c^q(t)} \right)\right] \, \frac {dx}{\lambda_c(t)}.
\ee

Then,

\be
\begin{aligned}
    M_3^{(1)}
     & \le \frac{(\lambda_c'(t))^2}{2\mu_c(t)} \int \left(\frac x {\lambda_c(t)} \right)^2 \left|\phi' \left( \frac x {\lambda_c(t)} \right)\right| \, \frac {dx}{\lambda_c(t)} \\
     & = \frac{(\lambda_c'(t))^2}{2\mu_c(t)} \int y^2 \left| \phi' (y)\right| \, dy                                                                                       \\
     & = C \frac{(\lambda_c'(t))^2}{2\mu_c(t)} \lesssim C \frac{1}{t^{3-3b}\log^4t},
\end{aligned}
\ee

which is integrable for $t\gg1$ when $b < 2/3$.

For \eqref{M3}, note that $|\phi(x)| \le 1$, and $\int u^2 \, dx \le \norm{u}_{H^1}$. This gives us

\be
\begin{aligned}
    |M_3| & \le \left|\frac q {2\mu_c(t)} \frac {\lambda_c'(t)} {t^{1-c}\lambda_c(t)} \right| \int \left(\frac x {\lambda_c^q(t)} \right)^2 \left[ \phi'' \left( \frac x {\lambda_c^q(t)} \right)\right]^2  \, dx \\ &+ \left|\frac {qt^{1-c}} {2\mu_c(t)} \frac {\lambda_c'(t)} {\lambda_c(t)} \right|\int   u^2 \, dx \\
          & \le C\left|\frac {q\lambda_c^q(t)} {2\mu_c(t)} \frac {\lambda_c'(t)} {t^{1-c}\lambda_c(t)} \right| + C\left|\frac {qt^{1-c}} {2\mu_c(t)} \frac {\lambda_c'(t)} {\lambda_c(t)}\right|                      \\
          & \lesssim C \frac{1}{t^{3 - b(2+q)}\log^{2+q}t} +C \frac 1 {t\log^2 t},
\end{aligned}
\ee

which is integrable for $t \gg 1$ because $b \le \frac 2 {2+q}$.

If we combine our estimates for $M_1$, $M_2$, $M_3$, and $M_5$ with \eqref{M-deriv}, we get

\be
\begin{aligned}
    \frac d {dt} \mathcal M_q(t) + M_1 + &\frac C {t^{3-2b} \log^3 t} + M_3 + M_5 \\& \ge -\left|\frac{ 1}{2\mu_c(t)\lambda_c(t)}\right| \int \left[ \phi' \left( \frac x {\lambda_c(t)} \right)\phi' \left( \frac x {\lambda_c^q(t)} \right)\right] u^2 \, dx  + M_4 \\
        & =-N + M_4
\end{aligned}
\ee

If we assume $t \gg 1$, we get

\be
\begin{aligned}
    -N + M_4 & = -\left|\frac{ 1}{2\mu_c(t)\lambda_c(t)}\right| \int \left[ \phi' \left( \frac x {\lambda_c(t)} \right)\phi' \left( \frac x {\lambda_c^q(t)} \right)\right] u^2 \, dx                                                   \\&+ \frac 1 {\mu_c(t)\lambda_c(t)} \int \phi' \left( \frac x {\lambda_c(t)} \right) \phi' \left( \frac x {\lambda_c^q(t)} \right) \frac {u^2} 2 \, dx \\
    &+ \frac 1 {\mu_c(t)\lambda_c(t)} \int \phi' \left( \frac x {\lambda_c(t)} \right) \phi' \left( \frac x {\lambda_c^q(t)} \right) \left[G \left( \frac b 2 u^2 + \frac {3-b} 2 (\p_x u)^2\right)\right] \, dx\\
             & = \underbrace{ \frac 1 {\mu_c(t)\lambda_c(t)} \int \phi' \left( \frac x {\lambda_c(t)} \right) \phi' \left( \frac x {\lambda_c^q(t)} \right) \left[G \left( \frac b 2 u^2 + \frac {3-b} 2 (\p_x u)^2\right)\right] \, dx}_{(i)} \\
             & =  \underbrace{ \frac 1 {\mu_c(t)\lambda_c(t)} \int G \left[ \phi' \left( \frac x {\lambda_c(t)} \right) \phi' \left( \frac x {\lambda_c^q(t)} \right)\right] \left( \frac b 2 u^2 + \frac {3-b} 2 (\p_x u)^2\right)\, dx}_{(ii)} \\
             & \approx \frac 1 {t \log t} \int \phi' \left( \frac x {\lambda_c(t)} \right)\phi' \left( \frac x {\lambda_c^q(t)} \right)\left( \frac b 2 u^2 + \frac {3-b} 2(\p_x u)^2\right)\, dx
\end{aligned}
\ee

which is the quantity we need to estimate for our lemma. (i) and (ii) are equal because of Fubini's theorem since $Gf = \frac 1 2 e^{-|x|} * f$, $u^2 + (\p_x u)^2$ is $L^1(\R)$, and $ \phi' \left( \frac x {\lambda_c(t)} \right)\phi' \left( \frac x {\lambda_c^q(t)} \right)$ is a Schwartz function. For the last step, we used the fact that

\be
\phi' \left( \frac x {\lambda_c(t)} \right)\phi' \left( \frac x {\lambda_c^q(t)} \right) = e^{-\left| \frac x {\lambda} \right|} e^{- \left| \frac x {\lambda^q} \right|} = e^{-\left| x \frac {\lambda^{q-1}+1}{\lambda^q} \right|}, \quad \frac {\lambda^{q-1}+1}{\lambda^q} \approx \frac 1 \lambda
\ee

and Lemma \eqref{exp-est}.

Now if we look back at (3.6), we showed that

\be
\begin{aligned}
    \frac d {dt} \mathcal M_q(t) + M_1 + C \frac 1 {t^{3-2b} \log^3 t} + M_3 + M_5 \ge \\ \frac C {t \log t} \int \phi' \left( \frac x {\lambda_c(t)} \right)\phi' \left( \frac x {\lambda_c^q(t)} \right)\left( \frac b 2 u^2 + \frac {3-b} 2(\p_x u)^2\right)\, dx
\end{aligned}
\ee

The left hand side of the above inequality is integrable for $t \gg 1$. Since the right hand side is nonnegative, the right hand side of our inequality must be integrable when $t \gg 1$ as well.

Note that $\frac b 2 u^2 + \frac {3-b} 2 (\p_x u)^2 \ge \min\{b, \frac {3-b} 2\} (u^2 + (\p_x u))^2$, so we know that
\[\frac C {t \log t} \int \phi' \left( \frac x {\lambda_c(t)} \right)\phi' \left( \frac x {\lambda_c^q(t)} \right)\left( u^2 + (\p_x u)^2\right)\, dx\]
is integrable for $t \gg 1$, which proves the lemma.

\ep

Now we can prove the weak decay estimate

\subsection{Proof of Theorem \ref{TH3}}
Lemma 3.1 states that

\be
\int_{10}^\infty \frac 1 {t \log t} \int \phi' \left( \frac x {\lambda_c(t)} \right)\phi' \left( \frac x {\lambda_c^q(t)} \right) (u^2 + (\p_x u)^2) \, dxdt < +\infty
\ee

However, $t \ln t$ is not integrable for any $t > 1$. This implies that

\be
\liminf_{t \to +\infty} \int \phi' \left( \frac x {\lambda_c(t)} \right)\phi' \left( \frac x {\lambda_c^q(t)} \right)  (u^2 + (\p_x u)^2) \, dx = 0
\ee

Recall that $\lambda_c(t) = \frac {t^c}{\log t}$. This means that when $|x| \in (-\Lambda_c(t), \Lambda_c(t))$, there exists $\delta > 0$ such that $\phi' \left( \frac x {\lambda_c(t)} \right)\phi' \left( \frac x {\lambda_c^q(t)} \right) \ge \delta$. This then implies that

\be
\liminf_{t \to +\infty}\norm{u}_{H^1(\Lambda_c(t))} = 0
\ee

which gives us the weak decay estimate.

\section{Proof of Theorem \ref{TH1}}

\subsection{First Functional Estimate}

Consider the functional

\be
\label{Iv}
\mathcal I(t) = \int \phi\left( \frac x {\lambda_c(t)} \right) u(t,x) \, dx
\ee

where $\phi(x) = \tanh x$. Note that \eqref{conserved} implies that $\mathcal I(t)$ is bounded in time. Now we wish to prove

\begin{lem}
    \label{weak-decay}
    \be
    \label{norm-decay}
    \int_t^\infty \frac 1 {\lambda_c(t)} \int \phi'\left( \frac x {\lambda_c(t)} \right) (u^2 + (\p_x u)^2) \, dx < \infty,
    \ee

    for $t \gg 1$.
\end{lem}

\bp
First, take the time derivative of \eqref{Iv}. Applying \eqref{bfam} and integrating by parts gives us

\be
\label{Ivderiv}
\begin{aligned}
    \frac d {dt} \mathcal I(t) &= - \frac {\lambda_c'(t)} {\lambda_c(t)} \int \frac x {\lambda_c(t)} \phi' \left( \frac x {\lambda_c(t)} \right) u \, dx + \int \phi \left( \frac x {\lambda_c(t)} \right) u_t \, dx \\
    &= - \frac {\lambda_c'(t)} {\lambda_c(t)} \int \frac x {\lambda_c(t)} \phi' \left( \frac x {\lambda_c(t)} \right) u \, dx - \int \phi \left( \frac x {\lambda_c(t)} \right) \p_x \left( \frac {u^2} 2 \right) \, dx \\
    & - \int \phi \left( \frac x {\lambda_c(t)} \right) \p_x \left( G \left( \frac b 2 u^2 + \frac {3-b} 2 (\p_x u)^2 \right) \right) \, dx \\
    & = - \frac {\lambda_c'(t)} {\lambda_c(t)} \int \frac x {\lambda_c(t)} \phi' \left( \frac x {\lambda_c(t)} \right) u \, dx + \frac 1 {2\lambda_c(t)} \int \phi' \left( \frac x {\lambda_c(t)} \right) u^2 \, dx \\
    & + \frac x {\lambda_c(t)} \int \phi' \left( \frac x {\lambda_c(t)} \right) G \left( \frac b 2 u^2 + \frac {3-b} 2 (\p_x u)^2 \right) \, dx.
\end{aligned}
\ee

If we let $c=1/2$, Young's inequality, a change of variables, and \eqref{lu-deriv} tells us that

\be
\label{Iv-log}
\begin{aligned}
    \left| \frac {\lambda_c'(t)} {\lambda_c(t)} \int \frac x {\lambda_c(t)} \phi' \left( \frac x {\lambda_c(t)} \right) u \, dx \right| &\le \frac {(\lambda_c'(t))^2} {2\delta} \int \left( \frac x {\lambda_c(t)} \right)^2 \phi' \left( \frac x {\lambda_c(t)} \right) \, \frac {dx} {\lambda_c(t)} \\
    & + \frac \delta {2\lambda_c(t)} \int \phi' \left(\frac x {\lambda_c(t)} \right) u^2 \, dx \\
    & \lesssim \frac C {\delta t \log^2 t} + \frac \delta {2\lambda_c(t)} \int \phi' \left(\frac x {\lambda_c(t)} \right) u^2 \, dx,
\end{aligned}
\ee

 where $\delta > 0$ is a constant that will be chosen later. Plugging this into \eqref{Ivderiv} gives us

\be
\begin{aligned}
    \frac d {dt} \mathcal I(t) & \gtrsim -\frac C {t \log^2 t} -\frac \delta {2\lambda_c(t)} \int \phi' \left(\frac x {\lambda_c(t)} \right) u^2 \, dx \\ 
    &+ \frac 1 {\lambda_c(t)} \underbrace{\int \phi' \left( \frac x {\lambda_c(t)} \right) G \left( \frac b 2 u^2 + \frac {3-b} 2 (\p_x u)^2 \right) \, dx }_{(i)}\\
    & = -\frac C {t \log^2 t} -\frac \delta {2\lambda_c(t)} \int \phi' \left(\frac x {\lambda_c(t)} \right) u^2 \, dx \\
    & + \frac 1 {\lambda_c(t)}\underbrace{ \int G\left[ \phi' \left( \frac x {\lambda_c(t)} \right) \right] \left( \frac b 2 u^2 + \frac {3-b} 2 (\p_x u)^2 \right) \, dx}_{(ii)} \\
    & \gtrsim -\frac C {t \log^2 t} -\frac \delta {2\lambda_c(t)} \int \phi' \left(\frac x {\lambda_c(t)} \right) u^2 \, dx \\
    & + \frac {C(\lambda)} {\lambda_c(t)} \int \phi' \left( \frac x {\lambda_c(t)} \right) \left( \frac b 2 u^2 + \frac {3-b} 2 (\p_x u)^2 \right) \, dx \\
    & \gtrsim \frac C {t \log^2 t} + \frac {1} {2\lambda_c(t)} \int \phi' \left( \frac x {\lambda_c(t)} \right) \left( \epsilon u^2 + \frac {3-b} 2 (\p_x u)^2 \right) \, dx, \, \epsilon > 0.
\end{aligned}
\ee

(i) and (ii) are equal because of Fubini's theorem since $Gf = \frac 1 2 e^{-|x|} * f$, $u^2 + (\p_x u)^2$ is $L^1(\R)$, and $ \phi' \left( \frac x {\lambda_c(t)} \right)$ is a positive Schwartz function. The second to last inequality comes from \eqref{exp-est} and Fubini's theorem to apply $G$ to $\phi'$ instead. Recall that $C(\lambda) = \frac {\lambda^2} {\lambda^2 - 1}$. The last inequality comes from an appropriate choice of $\delta$ as long as $\lambda$ is large enough. Since $\mathcal I(t)$ is bounded, and $\frac 1 {t \log^2 t}$ is integrable for $t \gg 1$, we get 

\be
\int_t^\infty \frac 1 {\lambda_c(t)} \int \phi'\left( \frac x {\lambda_c(t)} \right) \left(\epsilon u^2 + \frac {3-b} 2 (\p_x u)^2\right) \, dx < \infty,
\ee

where the integrand clearly is bounded below by a constant multiple of the left hand side of \eqref{norm-decay}. This proves the lemma. Note that $\frac 1 {\lambda_c(t)}$ is not integrable in time. This implies that 

\be
\liminf_{t \to \infty} \int \phi' \left( \frac x {\lambda_c(t)} \right) \left( u^2 + (\p_x^2 u)^2 \right) \, dx = 0.
\ee

Otherwise the previous integrand could not be integrable.

\ep

\subsection{Second Functional Estimate}

To complete the proof of Theorem \ref{TH1}, we will need to consider the following functional

\be
\label{EnergyFunctional}
    \mathcal E (t) \coloneqq \int \psi\left( \frac x {\lambda_c(t)} \right) m(t,x) \, dx,
\ee

where $\psi(x) = \mathrm{sech}^2 \, x$. $\mathcal E(t) < +\infty$ for all $t$ since $\int m \, dx$ is a constant of motion, and $\mathrm{sech}^2(x)$ is a Schwartz function. Note that $\psi'(x) = \mathrm{sech}\, x \tanh x$, so $|\psi ' (x) | \le \psi(x)$.

First we need to compute the derivative of the above functional. If we apply \eqref{cons-bfam}, we get

\be
\begin{aligned}
    \label{e-deriv}
    \frac d {dt} \mathcal E(t) &= -\frac{\lambda_c'(t)}{\lambda_c(t)} \int \frac x {\lambda_c(t)} \psi'\left(\frac x {\lambda_c(t)} \right) m \, dx + \int \psi \left( \frac x {\lambda_c(t)} \right) m_t \, dx \\
    &= -\frac{\lambda_c'(t)}{\lambda_c(t)} \int \frac x {\lambda_c(t)} \psi'\left(\frac x {\lambda_c(t)} \right) (u- \p_x^2 u) \, dx \\
    &- \int \psi \left(\frac x {\lambda_c(t)} \right)  (u \p_x u - u \p_x^2 u) \, dx \\
    &- b\int \psi \left(\frac x {\lambda_c(t)} \right) \left( u \p_x u - \p_x u \p_x^2 u \right) \, dx 
\end{aligned}
\ee

Now consider the following lemma

\begin{lem}
\label{e-bound}
\be
    \left| \frac d {dt} \mathcal E(t) \right| \lesssim \frac 1 {t\log^2t} + \frac {\log t} {t^{3/2}} + \frac 1 {\lambda_c(t)} \int \psi \left( \frac x {\lambda_c(t)} \right) (u^2 + (\p_x u) ^2) \, dx.
\ee
\end{lem}

\bp

If we look at \eqref{e-deriv}, we can integrate by parts to get
\be
\begin{aligned}
    \frac d {dt} \mathcal E(t) &= -\frac{\lambda_c'(t)}{\lambda_c(t)} \int \frac x {\lambda_c(t)} \psi'\left(\frac x {\lambda_c(t)} \right) (u- \p_x^2 u) \, dx \\
    &- \int \psi \left(\frac x {\lambda_c(t)} \right)  (u \p_x u - u \p_x^3 u) \, dx \\
    &- b\int \psi \left(\frac x {\lambda_c(t)} \right) \left( u \p_x u - \p_x u \p_x^2 u \right) \, dx \\
    &= -\frac{\lambda_c'(t)}{\lambda_c(t)} \int \frac x {\lambda_c(t)} \psi'\left(\frac x {\lambda_c(t)} \right) u \, dx - \frac{\lambda_c'(t)}{\lambda_c(t)} \int \frac x {\lambda_c(t)} \psi'\left(\frac x {\lambda_c(t)} \right)  \p_x^2 u \, dx \\
    &+ \frac 1 {2\lambda_c(t)} \int \psi' \left(\frac x {\lambda_c(t)} \right) (u^2 + (\p_x u)^2 ) \, dx - \frac 1 {\lambda_c(t)} \int  \psi' \left(\frac x {\lambda_c(t)} \right)  u \p_x^2 u \, dx \\
    &+ \frac b {2\lambda (t)} \int \psi' \left(\frac x {\lambda_c(t)} \right) \left( u^2 - (\p_x u)^2  \right) \, dx.
\end{aligned}
\ee

Now if we integrate by parts, apply Young's inequality, and use \eqref{lu-deriv} for $b= 1/2$, we can show that
\be
\label{m-youngs}
\begin{aligned}
    &\left| -\frac{\lambda_c'(t)}{\lambda_c(t)} \int \frac x {\lambda_c(t)} \psi'\left(\frac x {\lambda_c(t)} \right) (u- \p_x^2 u) \, dx \right| = \left| -\frac{\lambda_c'(t)}{\lambda_c(t)} \int \frac x {\lambda_c(t)} \psi'\left(\frac x {\lambda_c(t)} \right) u \, dx \right| \\
    &+ \left| \frac{\lambda_c'(t)}{(\lambda_c(t))^2} \int \left[ \frac x {\lambda_c(t)} \psi'' \left( \frac x {\lambda_c(t)} \right) + \psi'\left(\frac x {\lambda_c(t)} \right) \right] \p_x u \, dx \right| \\
    &\le \frac {(\lambda_c'(t))^2} {2} \int \left( \frac x {\lambda_c(t)} \right)^2 \psi \left( \frac x {\lambda_c(t)} \right) \, \frac {dx} {\lambda_c(t)}
    + \frac 1 {\lambda_c(t)} \int  \psi \left( \frac x {\lambda_c(t)} \right) u^2 \, dx + \frac {D \log t} {t^{3/2}} \\
    &\le \frac {C} {t \log^2 t}+ \frac 1 {\lambda_c(t)} \int \psi \left( \frac x {\lambda_c(t)} \right) u^2 \, dx + \frac {D \log t} {t^{3/2}}.
\end{aligned}
\ee

Recall that $|\psi'(x)| \le \psi(x)$. Now if we plug in \eqref{m-youngs} into \eqref{e-deriv} and apply the triangle inequality, we get 

\be
\begin{aligned}
    \left| \frac d {dt} \mathcal E(t) \right| &\le\frac {E} {t \log^2 t}+ \frac {D \log t} {t^{3/2}} + \frac C {\lambda_c(t)} \int \psi \left( \frac x {\lambda_c(t)} \right) (u^2 + (\p_x u)^2) \, dx \\
    & +  \frac 1 {\lambda_c(t)} \int  \psi \left(\frac x {\lambda_c(t)} \right)  u \p_x^2 u \, dx \\
    &\le \frac {1} {t \log^2 t}+ \frac {D \log t} {t^{3/2}} + \frac C {\lambda_c(t)} \int \psi \left( \frac x {\lambda_c(t)} \right) (u^2 + (\p_x u)^2) \, dx \\
    & +  \frac 1 {(\lambda_c(t))^2} \int  \psi \left(\frac x {\lambda_c(t)} \right) u \p_x u \, dx + \frac 1 {\lambda_c(t)} \int  \psi \left(\frac x {\lambda_c(t)} \right) (\p_x u)^2 \, dx \\
    &\lesssim \frac 1 {t\log^2t} + \frac {\log t} {t^{3/2}} + \frac 1 {\lambda_c(t)} \int \psi \left( \frac x {\lambda_c(t)} \right) (u^2 + (\p_x u) ^2) \, dx
\end{aligned}
\ee

where the last inequality simply comes from Young's ienquality.

\ep

Now we can prove Theorem \ref{TH1}

\subsection{Proof of Theorem \ref{TH1}}

\bp

Say $t_0 \gg 1$. The right hand side of \eqref{e-bound} is integrable based on Lemma \ref{weak-decay} when we choose $c=1/2$ for $\lambda_c(t) = \frac {t^c} {\log t}$ and integrate from $t_0$ to $\infty$. This tells us that

    \be
    \begin{aligned}
        |\mathcal E(t_0)| &= \left| L+ \int_{t_0}^\infty \frac d {dt} \mathcal E(t) \, dt \right| \le L + \int_{t_0}^\infty \left| \frac d {dt} \mathcal E(t) \right| \, dt \\ &\lesssim \int_{t_0}^{\infty} \frac 1 {t\log^2t} + \frac {\log t} {t^{3/2}} + \frac 1 {\lambda_c(t)} \int \psi \left( \frac x {\lambda_c(t)} \right) \left(u^2 + (\p_x u)^2 \right) \, dxdt,
    \end{aligned}
    \ee

which then implies there exists a limit for $\mathcal E(t)$ as $t \to \infty$ by Dominated Convergence since $\mathcal E(t) \ge 0$. However, note that we do not know whether or not the limit is 0.

Recall that $m \ge 0$. With the definition of $\mathcal E(t)$, we have

\be
\mathcal E(t) = \int \psi \left( \frac x {\lambda_c(t)} \right) m (t,x)\, dx \to L.
\ee

We also have $1 \ge \psi(x/\lambda_c(t)) = e^{-|x / \lambda_c(t)|} \ge e^{-1}$ for all $x \in (-\Lambda_{1/2}(t), \Lambda_{1/2}(t)$. This implies that

\be
    e\mathcal E(t) \ge \norm{m(t)}_{L^1(\Lambda_{1/2}(t))} \ge \mathcal E(t)
\ee.

Note that we essentially get the same estimate for $\norm{u}_{L^1(\Lambda_{1/2})}$. To see this, recall that $m = u - \p_x^2 u$. This tells us that
\be
\begin{aligned}
    \mathcal E(t) &=\int \psi \left( \frac x {\lambda_c(t)} \right) m \, dx \\
    & = \int \psi \left( \frac x {\lambda_c(t)} \right) u \, dx - \int \psi \left( \frac x {\lambda_c(t)} \right) \p_x^2 u \, dx \\
    &= \int \psi \left( \frac x {\lambda_c(t)} \right) u \, dx - \frac{1} {\lambda^2(t)}\int \psi'' \left( \frac x {\lambda_c(t)} \right) u \, dx.
\end{aligned}
\ee

Since $|\psi''(x)| \le \psi(x)$, $\lambda_c(t) \to \infty$ as $t \to \infty$, and $\mathcal I(t)$ is bounded, the second term in the last inequality goes to 0. $u$ is signed since $u = Gm$, so we get

\be
\begin{aligned}
L = \lim_{t \to \infty} \mathcal E(t) &= \lim_{t \to \infty} \int \psi \left( \frac x {\lambda_c(t)} \right) u \, dx - \lim_{t \to \infty}\frac{1} {\lambda^2(t)}\int \psi'' \left( \frac x {\lambda_c(t)} \right) u \, dx\\
 &= \lim_{t \to \infty} \int \psi \left( \frac x {\lambda(t)} \right) u \, dx.
\end{aligned}
\ee

By a similar argument as before, this implies that 
\be
    e\mathcal E(t) \ge \norm{u(t)}_{L^1(\Lambda_{1/2}(t))} \ge \mathcal E(t)
\ee,

which then implies that

\be
\lim_{t \to \infty} \norm{u(t)}_{L^1(\Lambda_{1/2}(t))} = \lim_{t \to \infty} \norm{m(t)}_{L^1(\Lambda_{1/2}(t))} = L.
\ee

It is difficult to get any further results from here without further assumptions. The proof in \cite{ACKM2019} makes use of the fact that there exists a subsequential limit for a similar functional. From there, integrability would imply that limit, 0, would be the only possible limit. However, we do not know that 0 is the limit. In fact, we can't even determine the limit for the $H^1$ norm on $\Lambda_{1/2}(t)$ because the $L^1$ norm of $m$ only gives an upper bound for the $H^1$ norm of $u$.

\ep

\section{Proof of Theorem \ref{TH4}}

\subsection{Proof of Functional Estimate}

We will follow the method of proof used in \cite{K2019}. Consider the following functional

\be
\label{travel-decay}
\mathcal I(t) = \int \phi \left( \frac{x-\sigma t} L \right) m (t,x)\, dx
\ee

where $\sigma, \, L > 0$. By our assumptions, this functional is bounded in time, and if we take $\phi(x) = \frac 1 2 (\tanh x + 1)$, our functional is nonnegative as well given $m$ is nonnegative. To prove the result, we want to show 

\begin{lem}
    \label{travel-est}
    There exists some $\sigma > 0$ such that.
    \be
    \frac d {dt} \mathcal I(t) \lesssim -\frac \sigma L \int \phi' \left( \frac {x - \sigma t} L \right) (u^2 + (\p_x u)^2) \, dx
    \ee
\end{lem}

\bp

Following a similar method to the proof of Lemma \ref{e-bound}, we first take the derivative of $\mathcal I(t)$, and get

\be
\begin{aligned}
    \frac d {dt} \mathcal I(t) &= - \frac \sigma L \int \phi' \left( \frac{x - \sigma t} L \right) m \, dx + \int \phi \left( \frac{x - \sigma t} L \right) m_t \, dx \\
    &= - \frac \sigma L \int \phi' \left( \frac{x - \sigma t} L \right) m \, dx - \int \phi \left( \frac{x - \sigma t} L \right) (u\p_x m + b \p_x u m) \, dx \\
    &= - \frac \sigma L \int \phi' \left( \frac{x - \sigma t} L \right) m \, dx - \int \phi \left( \frac{x - \sigma t} L \right) \left( u \p_x u - u \p_x^3 u \right) \, dx \\ 
    & - b\int \phi \left( \frac{x - \sigma t} L \right)  \left( u \p_x u - \p_x u \p_x^2 u \right) \, dx.
\end{aligned}
\ee

Then we integrate by parts to get

\be
\begin{aligned}
    \frac d {dt} \mathcal I(t) &= - \frac \sigma L \int \phi' \left( \frac{x - \sigma t} L \right) m \, dx - \frac \sigma {2L} \int \phi' \left( \frac{x - \sigma t} L \right) u^2 \, dx \\ 
    & + \frac \sigma {L} \int \phi' \left( \frac{x - \sigma t} L \right) u \p_x^2 u \, dx - \int \phi \left( \frac{x - \sigma t} L \right) \p_x u \p_x^2 u \, dx \\
    & - b\frac \sigma {2L} \int \phi \left( \frac{x - \sigma t} L \right) \left( u^2- (\p_x u)^2 \right) \, dx.
\end{aligned}
\ee

Further integration by parts and collecting like terms gives us

\be
\begin{aligned}
    \frac d {dt} \mathcal I(t) &= - \frac \sigma L \int \phi' \left( \frac{x - \sigma t} L \right) m \, dx - \frac \sigma {2L} \int \phi' \left( \frac{x - \sigma t} L \right) u^2  \, dx \\ 
    & - \frac 3 2 \frac \sigma {L} \int \phi' \left( \frac{x - \sigma t} L \right) (\p_x u)^2\, dx + \frac {\sigma^2} {L^2} \int \phi'' \left( \frac{x - \sigma t} L \right) u \p_x u \, dx \\
    & - b\frac \sigma {2L} \int \phi \left( \frac{x - \sigma t} L \right) \left( u^2- (\p_x u)^2 \right) \, dx \\
    & \lesssim - \frac \sigma L \int \phi' \left( \frac{x - \sigma t} L \right) m \, dx - \frac {C-b} 2 \frac {\sigma} L \int \phi' \left( \frac{x - \sigma t} L \right) (u^2 + (\p_x u)^2) \, dx \\
    &+ \frac 1 2 \frac {\sigma^3} {L^3} \int \phi''' \left( \frac {x - \sigma t} L \right) u^2 \, dx,
\end{aligned}
\ee

Note that $3 - b > 0$ for $b \in (0, 3)$. \eqref{p-estimate} implies we can also get the bound

\be
\frac d {dt} \mathcal I(t) \lesssim - \left(\frac \sigma L - \frac 1 2 \frac {\sigma^3} {L^3}\right) \int \phi' \left( \frac{x - \sigma t} L \right) m \, dx -\frac {C\sigma} L \int \phi' \left( \frac{x - \sigma t} L \right) (u^2 + (\p_x u)^2) \, dx
\ee

We can simply let $L$ be sufficiently larger than $\sigma$ to get the result. This guarantees that the first term in the above inequality is negative. Then we can simply ignore the first negative term to get the result.

\ep

\subsection{Proof of Theorem \ref{TH4}}

\bp

Consider the following functional

\be
I_{t_0} (t) = \int \phi \left( \frac {x - \sigma t_0 + \frac \sigma 2 (t_0 - t)} {L} \right) m \, dx,
\ee

where $\phi(x) = \frac 1 2 (\tanh x + 1)$, and $t_0 > 2$.

Lemma \ref{travel-est} implies that

\be
\frac d {dt} I_{t_o} (t) \lesssim - \int \phi' \left(\frac {x - \sigma t_0 + \frac \sigma 2 (t_0 - t)} {L} \right) (u^2 + (\p_x u)^2) \, dx.
\ee

Note that $\phi'(x) = \text{sech}^2 x$. Since the integrand is nonnegative, we get that $I_{t_0}(t)$ is decreasing on the interval $t \in [2, t_0]$. Moreover, since $\phi(x)$ tends to 0 as $x \to -\infty$, so $\phi\left( \frac {x-at-b} L \right) \rightharpoonup 0$ for positive $a, b$. This implies

\be
\limsup_{t \to \infty} \int \phi \left( \frac {x - a t - b} {L} \right) m (\delta, x) \, dx = 0
\ee

when $a, b, \delta > 0$. Since $I_{t_0}(t)$ is decreasing on $[2, t_0]$, we get the following

\be
0 \le \int \phi \left( \frac {x - \sigma t_0} {L} \right) m (t_0, x)\, dx \le \int \phi \left( \frac {x - \frac \sigma 2 t_0 - \sigma} {L} \right) m (2, x) \, dx \xrightarrow{t_0 \to \infty} 0.
\ee

Recall that $1 \ge \phi \left( \frac {x - \frac \sigma t_0} {L} \right) \ge 1/2$ for $x \ge \frac \sigma t_0$. Since $m(t, x) \ge 0$, we can apply \eqref{p-estimate} to get

\be
    \frac 1 2 \| u(t) \|_{H^1(\tilde{\sigma}, \infty)} \le \frac 1 2 \| m(t) \|_{L^1(\tilde{\sigma}t, \infty)} \le I_{t_0}(t) \xrightarrow{t \to \infty} 0
\ee

In fact, \eqref{p-estimate} implies we get similar asymptotic decay for $u \in W_1^p(\tilde{\sigma}, \infty)$. This completes the proof.

\ep

\bibliographystyle{abbrv}
\bibliography{ref}

\end{document}